\documentclass{amsart}

\usepackage{color}

\newtheorem{theorem}{Theorem}[section]
\newtheorem{lemma}[theorem]{Lemma}
\newtheorem{prop}[theorem]{Proposition}

\theoremstyle{definition}
\newtheorem{definition}[theorem]{Definition}
\newtheorem{example}[theorem]{Example}
\newtheorem{conjecture}[theorem]{Conjecture}
\newtheorem{cor}[theorem]{Corollary}

\theoremstyle{remark}

\numberwithin{equation}{section}

\newcommand{\bQ}{ {\mathbb Q}}
\newcommand{\bA}{ {\mathbb A}}

\newcommand{\bC}{ {\mathbb C}}

\newcommand{\bF}{ {\mathbb F}}
\newcommand{\bN}{ {\mathbb N}}
\newcommand{\bZ}{ {\mathbb Z}}

\newcommand{\cD}{ {\mathcal D}}

\newcommand{\vx}{ { \mathbf{x}}}

\newcommand{\hd}{ {\odot}}

\begin{document}

\title{Power Series with Coefficients from a Finite Set}

\author{Jason Bell}
\address{Department of Pure Mathematics, University of Waterloo,
Waterloo, Ontario N2L 3G1, Canada}
\email{jpbell@uwaterloo.ca}
\thanks{Jason P.\ Bell was supported by the NSERC Grant RGPIN-2016-03632.}

\author{Shaoshi Chen}
\address{KLMM, Academy of Mathematics and Systems Science, Chinese Academy of Sciences, Beijing, 100190, China;
and Symbolic Computation Group, David R.\ Cheriton School of Computer Science, University of Waterloo,
Waterloo, Ontario N2L3G1, Canada}
\email{schen@amss.ac.cn}
\thanks{
Shaoshi Chen was supported by the NSFC grant 11501552 and
by the President Fund of the Academy of
Mathematics and Systems Science, CAS (2014-cjrwlzx-chshsh). This work was
also supported by the Fields Institute's 2015 Thematic Program on Computer
Algebra.}


\date{}



\begin{abstract}
We prove in this paper that a multivariate D-finite power series with coefficients from a finite set is rational.
This generalizes a rationality theorem of van der Poorten and Shparlinski in 1996.
\end{abstract}

\maketitle


\section{Introduction } \label{SECT:intro}
In his thesis~\cite{Hadamard1892}, Hadamard began the study of the
relationship between the coefficients of a power series and the properties of the function it
represents, especially its singularities and natural boundaries. Two special cases of the problem
have been extensively studied: one is on power series with integer coefficients and the other is
on power series with finitely many distinct coefficients.

In the first case, Fatou~\cite{Fatou1906} in 1906 proved a lemma on rational power series with integer coefficients,
which is now known as Fatou's lemma~\cite[p.\ 275]{Stanley1986}.
The next celebrated result is the P\'{o}lya-Carlson theorem, which asserts that a power series with integer coefficients
and of radius of convergence 1 is either rational or has the unit circle as its natural boundary.
This theorem was first conjectured in 1915 by P\'{o}lya~\cite{Polya1916} and later proved in 1921 by Carlson~\cite{Carlson1921}.
Several extensions of the P\'{o}lya-Carlson theorem have been presented in~\cite{Polya1921, Petersson1931, Gerig1969, Seinov1971, Martineau1971, Straube1987, Bell2014}.

In the second case, Fatou \cite{Fatou1906} was also the first to investigate power series with coefficients from a finite set.
The study was continued by P\'olya~\cite{Polya1916} in 1916, Jentzsch~\cite{Jentzsch1917} in 1917, Carlson~\cite{Carlson1918} in 1918 and finally Szeg\H{o} \cite{Szego1922, SzegoPapers1982} in 1922 settled the question by proving the following beautiful theorem (see~\cite[Chap.\ 11]{Remmert1998} and~\cite[Chap.\ 10]{Dienes1957} for its proof and related results).

\begin{theorem}[Szeg\H{o}, 1922]
Let~$F = \sum f(n) x^n$ be a power series with coefficients from a finite values of~$\bC$. If $F$ is continuable
beyond the unit circle then it is a rational function of the form~$F = P(x)/(1-x^m)$,
where~$P$ is a polynomial and~$m$ a positive integer.
\end{theorem}
Szeg\H{o}'s theorem was generalized in 1945 by Duffin and Schaeffer~\cite{Duffin1945} by assuming a weaker condition
that~$f$ is bounded in a sector of the unit circle. In 2008, P.\ Borwein et al. in~\cite{Borwein2008} gave a shorter proof of Duffin
and Schaeffer's theorem. By using Szeg\H{o}'s theorem, van der Poorten and Shparlinski proved the following result~\cite{Poorten1996}.

\begin{theorem}[van der Poorten and Shparlinski, 1996]\label{THM:ps}
Let~$F = \sum f(n) x^n$ be a power series with coefficients from a finite values of~$\bQ$. If~$f(n)$ satisfies a linear recurrence equation
with polynomial coefficients, then $F$ is rational.
\end{theorem}
A univariate sequence~$f: \bN \rightarrow K$ is \emph{P-recursive} if it satisfies a linear recurrence equation
with polynomial coefficients in~$K[n]$. A power series~$F = \sum f(n)x^n$ is \emph{D-finite} if it satisfies
a linear differential equation with polynomial coefficients in~$K[x]$. By~\cite[Theorem 1.5]{Stanley1980}, a sequence
$a(n)$ is P-recursive if and only if the power series~$F := \sum f(n)x^n$ is D-finite. The notion of D-finite power series
can be generalized to the multivariate case~(see Definition~\ref{DEF:dfinite}). Our main result is the following
multivariate generalization of Theorem~\ref{THM:ps}.

\begin{theorem}\label{THM:main}
Let $K$ be a field of characteristic zero, and let $\Delta$ be a finite subset of $K$.  Suppose that $f:\mathbb{N}^d\to \Delta$ with $d\geq 1$ is such that $$F(x_1,\ldots ,x_d):=\sum_{(n_1,\ldots ,n_d)\in \mathbb{N}^d} f(n_1,\ldots ,n_d)x_1^{n_1}\cdots x_d^{n_d} \in K[[x_1,\ldots ,x_d]]$$ is $D$-finite.  Then $F$ is rational.
\end{theorem}
We note that a multivariate rational power series $$F(x_1, \cdots ,x_d) = \sum f(n_1,\ldots ,n_d)x_1^{n_1}\cdots x_d^{n_d}$$ with all coefficients in $\{0,1\}$ has a very restricted form. In particular, the set $E$ of $(n_1,\ldots ,n_d)\in \mathbb{N}^d$ for which $f(n_1,\ldots ,n_d)\neq 0$ is \emph{semilinear}; that is there exist $n\in \bN$ and finite subsets $V_0, \ldots V_n$ of~$\bN^d$, and $b_1, \ldots, b_n\in \bN^d$
such that
\begin{equation}\label{EQ:semilinear}
E = V_0 \bigcup \left\{ \bigcup_{i=1}^n \left(b_i + \sum_{v\in V_i} v\cdot \mathbb{N}\right)\right\}.
\end{equation}
Although this result is known, we are unaware of a reference and give a proof of this fact in Proposition \ref{LEM:semilinear}.

The remainder of this paper is organized as follows.
The basic properties of D-finite power series are
recalled in Section~\ref{SECT:dfinite}.
The proof of Theorem~\ref{THM:main} is given in Section~\ref{SECT:proof}.
In Section~\ref{SECT:gf}, we
present several applications of our main theorem on generating functions
over nonnegative points on algebraic varieties.


\section{D-finite power series} \label{SECT:dfinite}
Throughout this paper, we let~$\bN$ denote the set of all nonnegative integers.
Let~$K$ be a field of characteristic zero and let $K(\vx)$ be
the field of rational functions in several variables $\vx= x_1, \ldots, x_d$ over~$K$.
By $K[[\vx]]$ we denote the ring of formal power series in~$\vx$ over~$K$
and by $K((\vx))$ we denote the field of fractions of $K[[\vx]]$.
For two power series $F = \sum f(i_1, \ldots, i_d)x_1^{i_1} \cdots x_d^{i_d}$
and $G = \sum g(i_1, \ldots, i_d)x_1^{i_1} \cdots x_d^{i_d}$, the Hadamard product of~$F$ and~$G$ is defined by
\[F\hd G =  \sum f(i_1, \ldots, i_d)g(i_1, \ldots, i_d) x_1^{i_1} \cdots x_d^{i_d}.\]
Let~$D_{x_1}, \ldots, D_{x_d}$ denote the derivations on~$K((\vx))$ with respect to $x_1, \ldots, x_d$,
respectively.
\begin{definition}[\cite{Lipshitz1988}]\label{DEF:dfinite}
A formal power series $F(x_1, \ldots, x_d) \in K[[\vx]]$ is said to be~\emph{D-finite} over~$K(\vx)$
if the set of all derivatives $D_{x_1}^{i_1}\cdots D_{x_d}^{i_d}(F)$ with~$i_j \in \bN$ span a
finite-dimensional $K(\vx)$-vector subspace of $K((\vx))$. Equivalently, for each $i\in \{1, \ldots, d\}$, $F$ satisfies
a nontrivial linear partial differential equation of the form
\[\left \{p_{i, m_i} D_{x_i}^{m_i} + p_{i, m_1-1}D_{x_i}^{m_i-1} + \cdots + p_{i, 0}\right\} F = 0 \quad \text{with~$p_{i, j}\in K[\vx]$}.\]
\end{definition}
The notion of D-finite power series was first introduced in 1980 by Stanley \cite{Stanley1980}, and has since become ubiquitous
in algebraic combinatorics as an important part of the study of generating functions (see~\cite[Chap. 6]{Stanley1999}).  We recall
some closure properties of this class of power series.

\begin{prop}[\cite{Lipshitz1989}]\label{PROP:dfinite}
Let~$\cD$ denote the set of all D-finite power series in~$K[[\vx]]$.
Then
\begin{itemize}
\item[(i)] $\cD$ forms a subalgebra of~$K[[\vx]]$, i.e., if $F, G \in \cD$ and~$\alpha, \beta \in K$,
then $\alpha F + \beta G \in \cD$
and~$FG\in \cD$.
\item[(ii)] $\cD$ is closed under the Hadamard product, i.e., if $F, G \in \cD$, then~$F\hd G \in \cD$.
\item[(iii)] If~$F(x_1, \ldots, x_d)$ is D-finite, and $$\alpha_1(y_1, \ldots, y_d),\ldots ,\alpha_d(y_1,\ldots ,y_d)\in K[[y_1,\ldots ,y_d]]$$
are algebraic over~$K(y_1, \ldots, y_d)$ and the substitution makes sense, then $F(\alpha_1, \ldots, \alpha_d)$
is also D-finite over $K(y_1, \ldots, y_d)$. \\ In particular, if $F(x_1, \ldots, x_d)$ is D-finite and the
evaluation of $F$ at $x_d=1$ makes sense, then $F(x_1, \ldots, x_{d-1}, 1)$ is D-finite.
\end{itemize}
\end{prop}

The coefficients of a D-finite power series are highly structured. In the univariate case, a power series $f = \sum a(n)x^n$ is
D-finite if and only if the sequence $a(n)$ is P-recursive, i.e., it satisfies a linear recurrence equation with polynomial coefficients in~$n$~\cite{Stanley1980}.
The structure in the multivariate case is much more profound, which was explored by Lipshitz in~\cite{Lipshitz1989}.
We continue this exploration to study the position of nonzero coefficients. To this end, we recall a notion of size
in the semigroup $(\bN, +)$, introduced by Bergelson et al. \cite{Bergelson1998} for general semigroups.  A subset $S\subseteq \bN$ is \emph{syndetic} if there is some positive integer $C$ such that if $n\in S$ then $n+i\in S$ for some $i\in \{1,\ldots ,C\}$.
\begin{lemma} \label{LEM:syndetic}
Let $K$ be a field of characteristic zero and let
\[ G(x_1,\ldots ,x_d)=\sum_{(n_1,\ldots ,n_d)\in \mathbb{N}^d} g(n_1,\ldots ,n_d)x_1^{n_1}\cdots x_d^{n_d}\in K[[\vx]]\]
be a $D$-finite power series over~$K(\vx)$.  Then the set
$$\{n\in \mathbb{N} \mid \exists (n_1,\ldots ,n_{d-1})\in \mathbb{N}^{d-1}~{\rm such~that}~g(n_1,\ldots ,n_{d-1},n)\neq 0\}$$ is either finite or syndetic.
\end{lemma}
\begin{proof}
We let $L$ denote the field of fractions of $K[[x_1,\ldots ,x_{d-1}]]$. Then we may regard $G$
as a power series in $L[[x_d]]$ and it is $D$-finite in~$x_d$ over~$L(x_d)$ and it is straightforward
to see that the lemma reduces to the univariate case.
Thus we now assume that $G(x)=\sum g(n)x^n \in L[[x]]$ is $D$-finite.
Then there exist $m\ge 1$, distinct nonnegative integers $a_1=0, \ldots ,a_m$, and nonzero polynomials $P_1,\ldots ,P_m\in L[z]$ such that
\[\sum_{j=1}^m P_j(n) g(n+a_j)=0\]
for all sufficiently large $n$.
Then there is some $M$ such that $P_1(n)\cdots P_m(n)\neq 0$ for $n>M$.
If $m=1$ then we then see that $g(n)=0$ for $n>M$.  Thus we assume that $m>1$. Then if $n>M$ and $g(n)$ is nonzero then $g(n+a_j)$ is nonzero for some $1<j\leq m$ and so we see that the set of $n$ for which $g(n)$ is nonzero is syndetic.
\end{proof}

\section{Proof of the main theorem}\label{SECT:proof}
The proof of Theorem~\ref{THM:ps} by van der Poorten and Shparlinski is
based on the fact that any univariate D-finite power series represents
an analytic function with only finitely many poles~\cite{Stanley1980}, so it is impossible to
have the unit circle as its natural boundary. Then their result follows from Szeg\H{o}'s theorem.
The singularities of analytic functions represented by multivariate
D-finite power series are much more involved. It is not known how to
extend Szeg\H{o}'s theorem to the multivariate case. Thus new ideas are needed
in order to generalize Theorem~\ref{THM:ps} to the multivariate case.

Before the proof of our main theorem, we first prove a lemma about finitely generated $\bZ$-algebras.
\begin{lemma}\label{lem:finite}
Let $R$ be a finitely generated $\mathbb{Z}$-algebra that is an integral domain of characteristic zero and suppose that $x$ is a nonzero element of $R$.  Then there exists a finite set of prime numbers $\{p_1,\ldots ,p_m\}$ such that if $n$ is a positive integer such that $x\in nR$ then $n$ is an element of the multiplicative semigroup of $\mathbb{N}$ generated by $p_1,\ldots ,p_m$.
\end{lemma}
\begin{proof}
Let $U$ denote the group of units of $R$.  By a result of Roquette~\cite{Roquette1957}
(or see~\cite[page 39, Corollary]{Lang1962}) we have that $U$ is a finitely generated
abelian group and so $U_0$, the subgroup of $U$ generated by the rational numbers in $U$
is a finitely generated subgroup of $\mathbb{Q}^*$.  In particular, there exist prime
numbers $q_1,\ldots , q_t$ such that every positive rational number in $U$ is in the
multiplicative subgroup of $\mathbb{Q}^*$ generated by $q_1,\ldots ,q_t$.  Thus if $x$ is a unit and $x\in nR$ then $n$ is an integer unit of $R$ and hence in the semigroup generated by $q_1,\ldots ,q_t$.  Hence we may assume that $x$ is not a unit.

Let $Q_1\cap \cdots \cap Q_d$ be the primary decomposition of the ideal $xR$ and
let $P_j$ denote the radical of $Q_j$ for $j=1,\ldots ,d$.  Since each $P_j$ is (proper) prime ideal,
it follows that there is at most one prime number in $P_j$ for $j=1,\ldots ,d$.
We now let $\{p_1,\ldots ,p_m\}$ denote the union of the prime numbers in $P_1,\ldots ,P_d$ and $\{q_1,\ldots ,q_t\}$.

Now suppose that $x\in pR$ with $p$ a prime number.  We claim that $p$ must be in the finite set $\{p_1,\ldots ,p_m\}$.
To see this, we have $x=pr$ for some $r\in R$ and so for each $i\in \{1,\ldots ,d\}$ we have $pr\in Q_i$.
Since $Q_i$ is primary that gives that either $r\in Q_i$ or there is some $k\ge 1$ such that $p^k\in Q_i$.
Now if $r\in Q_i$ for every $i\in \{1,\ldots ,d\}$ then we see that $r$ is in the intersection of the $Q_i$
and hence in $xR$.  This then gives that $r=xa$ for some $a\in R$ and so we have
$$x=pr=pax.$$  Since $x$ is nonzero and $R$ is an integral domain, this then gives that $p\in U$ and so
by the above remarks we must have $p\in \{q_1,\ldots ,q_t\}\subseteq \{p_1,\ldots ,p_m\}$.

Alternatively, we have that there is some $i$ such that $p^k\in Q_i$ for some $k\ge 1$.
Then $p\in P_i$ and so $p\in \{p_1,\ldots ,p_m\}$.
To finish the proof, observe that if $n$ is not in the multiplicative semigroup generated
by $p_1,\ldots ,p_m$ then there is some prime $q\not\in \{p_1,\ldots ,p_m\}$ such that $q$ divides $n$ and so
$x\not \in nR$ since $nR \subseteq qR$ and $x\not\in qR$. The result follows.
\end{proof}

\begin{proof}[{\bf Proof of Theorem \ref{THM:main}}] We prove this by induction on $d$.
When $d=0$, $F$ is constant and there is nothing to prove.  We now suppose that the
result holds whenever $d<k$ and we consider the case when $d=k$.  Since $F$ is $D$-finite,
we have that $F(x_1,\ldots ,x_k)$ satisfies a nontrivial linear differential equation of the form
\[\sum_{j=0}^{\ell} P_j(x_1,\ldots ,x_k) \partial_{x_k}^j F = 0,\]
where $P_0,\ldots ,P_{\ell}$ are polynomials in~$K[x_1, \ldots, x_k]$.
Translating this into a relation for the coefficients of $F$, we see that there exists
some positive integer $N$ and polynomials
$Q_{a_1,\ldots ,a_k}(t)\in K[t]$ for $(a_1,\ldots ,a_k)\in \{-N,\ldots ,N\}^k$, not all zero, such that
\begin{equation}
\label{eq: Q}
\sum_{-N\le a_1,\ldots  ,a_k\le N} Q_{a_1,\ldots ,a_k}(n_k) f(n_1-a_1,\ldots ,n_k-a_k)=0
\end{equation} for all $(n_1,\ldots ,n_k)\in \mathbb{N}^k$, where we take $f(i_1,\ldots ,i_k)=0$
if some $i_j$ is negative.
By dividing our polynomials $Q_{a_1,\ldots ,a_k}(t)$ by $t^a$ for some nonnegative integer $a$
if necessary, we may assume that $q(a_1,\ldots ,a_k):=Q_{a_1,\ldots ,a_k}(0)$ is nonzero for
some $(a_1,\ldots ,a_k)\in \{-N,\ldots ,N\}^k$.
We now let $R$ denote the $\mathbb{Z}$-subalgebra of $K$
generated by $\Delta$ and by the coefficients of $Q_{a_1,\ldots ,a_k}(t)\in K[t]$
with $(a_1,\ldots ,a_k)\in \{-N,\ldots ,N\}^k$.  Then $R$ is finitely generated.
By construction, we have
\[\sum_{-N\le a_1,\ldots  ,a_k\le N} q(a_1,\ldots ,a_k) f(n_1-a_1,\ldots ,n_k-a_k)\in n_kR\]
for all $(n_1,\ldots ,n_k)\in \mathbb{N}^k$.
Now let $\Gamma$ denote the set of all numbers of the form
\[\sum_{-N\le a_1,\ldots  ,a_k\le N} q(a_1,\ldots ,a_k) s(a_1,\ldots ,a_k)\]
with $s(a_1,\ldots ,a_k)\in \Delta\cup\{0\}$.  Then $\Gamma$ is a finite set.
By Lemma \ref{lem:finite}, there is a finite set of prime numbers
$p_1,\ldots ,p_m$ such that for each nonzero $x\in \Gamma$ we have that
if $n$ is a positive integer with $x\in nR$ then $n$ is in the semigroup generated by $p_1,\ldots ,p_m$.  In particular,
\[\sum_{-N\le a_1,\ldots  ,a_k\le N} q(a_1,\ldots ,a_k) f(n_1-a_1,\ldots ,n_k-a_k)=0\]
whenever $n_k$ is not in the multiplicative semigroup generated by $p_1,\ldots ,p_m$. Equivalently,
\[G(x_1,\ldots ,x_k):=F(x_1,\ldots, x_k)  \left( \sum_{0\le a_1,\ldots  ,a_k\le N} q(a_1,\ldots ,a_k) x_1^{a_1}\ldots x_k^{a_k}\right)x_1^N\cdots x_k^N\]
has the property that $g(n_1,\ldots ,n_k)=0$ whenever $n_k\ge N$ and $n_k-N$ is not in the semigroup generated by $p_1,\ldots ,p_m$,
where $g(n_1,\ldots ,n_k)$ denotes the coefficient of $x_1^{n_1}\cdots x_k^{n_k}$ in $G(x_1,\ldots ,x_k)$.
Since $G$ is just $F$ multiplied by a polynomial, $G(x_1,\ldots ,x_k)$ is $D$-finite by
Proposition~\ref{PROP:dfinite}~(i); moreover, all coefficients of $G$ lie in the finite set $\Gamma$.
But now Lemma \ref{LEM:syndetic} gives that there is some positive integer $M$
such that $g(n_1,\ldots ,n_k)=0$ whenever $n_k>M$ since a translate of the multiplicative semigroup generated by $p_1,\ldots ,p_m$
cannot be syndetic.  Thus we have
\[G=\sum_{i=0}^M G_i(x_1,\ldots ,x_{k-1})x_k^i\]
for some power series $G_0,\ldots ,G_M\in K[[x_1,\ldots ,x_{k-1}]]$.  Then for $i\in \{0,\ldots ,M\}$,
we have that $G_i x_k^i$ is the Hadamard product of $G$ with
$x_k^i \prod_{j=1}^{k-1} (1-x_j)^{-1}$ and so each $G_i x_k^i$ is $D$-finite by Proposition~\ref{PROP:dfinite}~(ii).
Then specializing $x_k=1$ gives each $G_i$ is $D$-finite~by Proposition~\ref{PROP:dfinite}~(iii). Since each $G_i$ has coefficients in a finite set, we see by the induction hypothesis that each $G_i$ is rational and so $G$ is rational.  But this now gives that $F$ is rational by our definition of $G$, completing the proof.
\end{proof}

\section{Generating functions over nonnegative integer points on algebraic varieties}\label{SECT:gf}
Let~$V\subseteq \bA_{K}^d$ be an affine algebraic variety over an algebraically closed field~$K$ of characteristic zero.
We define the \emph{generating function over nonnegative integer points} on~$V$ by
\[ F_V(x_1, \ldots, x_d) := \sum_{(n_1, \ldots, n_d)\in V \cap \bN^d}x_1^{n_1} \cdots x_d^{n_d}.\]
Then one can ask the following questions about the properties of $F_V$ that often reflect the global geometric structure of~$V$:
\begin{enumerate}
\item  When $F_V$ is zero? This is Hilbert Tenth Problem when $K$ is the field of rational numbers. In 1970, Matiyasevich \cite{Matijasevich1970, Davis1973} proved that this problem is undecidable.

    \smallskip

\item  When $F_V$ is a polynomial? If so, $V$ has only finitely many nonnegative integer points.
Siegel's theorem on integral points answers this question for a smooth algebraic curve C of genus $g\geq 1$ defined over a number field $K$~\cite[Chap.\ 7]{Bombieri2006}.

    \smallskip

\item When $F_V$ is a rational function? This is always true when the variety $V$ is defined by linear polynomials with integer coefficients~\cite[Chap.\ 4]{Stanley1986}.

        \smallskip

\item  When $F_V$ is $D$-finite?  By our main theorem, we see that this question is the same as question (3), by
taking $f(n_1, \ldots, n_d) =1$ if~$(n_1, \ldots, n_d) \in V\cap \bN^d$ and~$f(n_1, \ldots, n_d)=0$ otherwise (see Corollary \ref{COR:algvar}).

   \smallskip
\item When $F_V$ satisfies an algebraic differential equation?More precisely, we say that a power series
$F(x_1, \ldots, x_d)\in K[[x_1, \ldots, x_d]]$ is \emph{differentially algebraic} if the transcendence
degree of the field generated by all of the derivatives~$D_{x_1}^{i_1}\cdots D_{x_d}^{i_d}(F)$ with~$i_j \in \bN$
over~$K(x_1, \ldots, x_d)$ is finite. If a power series is not differentially algebraic, then it is called~\emph{transcendentally transcendental}.
For a nice survey on transcendentally transcendental functions, see Rubel~\cite{Rubel1989}.
\end{enumerate}

\begin{cor}\label{COR:algvar}
Let $V\subseteq \mathbb{A}^d_K$ be an affine variety over an algebraically closed field $K$ of characteristic zero.  Then
the power series
\[ F_V(x_1, \ldots, x_d) := \sum_{(n_1,\ldots ,n_d)\in V\cap \mathbb{N}^d} x_1^{n_1}\cdots x_d^{n_d}\]
is $D$-finite if and only if it is rational.
\end{cor}

To show an application of this corollary, let us consider the linear system $A\vx = 0$, where $A$ is a $d\times m$ matrix with integer entries.
Let~$E$ be the set of all vectors~$(n_1, \ldots, n_d) \in \bN^d$ such that~$A\vx =0$. We now give
a proof of the following classical theorem in enumerative combinatorics.

\begin{theorem}[Theorem 4.6.11 in~\cite{Stanley1986}]\label{THM:rational}
The generating function
\[F_E(x_1, \ldots, x_d) := \sum_{{(n_1, \ldots, n_d)\in E}} x_1^{n_1}\cdots x_d^{n_d}\]
represents a rational function of $x_1, \ldots, x_d$.
\end{theorem}
\begin{proof}
By Corollary~\ref{COR:algvar}, it suffices to show that $F_E$ is D-finite.
We first recall a fact proved by Lipshitz in~\cite[p.\ 377]{Lipshitz1988}
that if the power series $G(\vx) = \sum g(n_1, \ldots, n_d)x_1^{n_1} \cdots x_d^{n_d}$ is D-finite and $C\subseteq \bN^d$
is the set of elements of $\mathbb{N}^d$ satisfying a finite set of inequalities of the form $\sum a_i n_i + b\geq 0$, where the $a_i, b\in \bZ$, then
the power series
\[H(\vx) := \sum_{(n_1, \ldots n_k)\in C} g(n_1, \dots, n_d)x_1^{n_1}\cdots x_d^{n_d}\]
is D-finite. Note that $R(x_1, \ldots, x_d) := \sum x_1^{n_1}\cdots x_d^{n_d} = 1/\prod_{i=1}^d(1-x_i)$ is
D-finite and any equality~$\sum a_i n_i=0$ is equivalent to two inequalities~$\sum a_i n_i \geq 0$ and~$\sum (-a_i) n_i \geq 0$. Then the D-finiteness of~$F_E$ follows from the fact.
\end{proof}

We now derive some properties of an algebraic variety~$E$ from the generating function~$F_E$ when~$d=2$.  We first prove a basic result that is probably well-known, but for which we are unaware of a reference.

\begin{prop}\label{LEM:semilinear}
Let~$$F(x_1, \cdots ,x_d) = \sum f(n_1,\ldots ,n_d)x_1^{n_1}\cdots x_d^{n_d}\in \bQ[[x_1,\ldots ,x_d]]$$ with~$f(n_1,\ldots ,n_d)\in \{0, 1\}$ for all $(n_1,\ldots ,n_d)\in \mathbb{N}^d$.
Then~$F$ is rational if and only if the support set~$E := \{(n_1,\ldots ,n_d)\in \bN^d\mid f(n_1, \ldots ,n_d) \neq 0\}$
of~$F$ is \emph{semilinear} (see Equation (\ref{EQ:semilinear}) for the definition of semilinearity).
\end{prop}
\begin{proof}
The sufficiency follows from Theorem~\ref{THM:rational}. For the other direction assume that $F(x_1, \ldots ,x_d)$
is rational. Then the sequence $f: \bN^d \rightarrow \{0, 1\}$ has a
rational generating function over any finite field $\bF_p$, where $p$ is a prime number.
(This follows from the fact that we can express the generating series for $F$ of the form $P/Q$ with $P$ and $Q$ polynomials in which the gcd of the collection of coefficients of $P$ and $Q$ is $1$.)
By Salon's theorem~\cite{Salon1987}, which is an multi-dimensional extension of the theorem by
Christol, Kamae, Mend\`es France, and Rauzy~\cite{Christol1980}, the sequence $f: \bN^d \rightarrow \{0, 1\}$ is $p$-automatic
for every prime number $p$. Then the Cobham-Semenov theorem~\cite{Durand2008}
implies that the support set $E$ of~$f$ is semilinear.
\end{proof}
We now use this result in the special case when $d=2$.
\begin{theorem} \label{THM:d=2}
Let~$p(x, y)\in K[x, y]$ be a nonzero polynomial satisfying that the generating function
\[F_p(x, y) := \sum_{\substack{(n, m)\in \bN^2 \\ p(n, m)=0}} x^ny^m\]
is rational.
Then $p = c\cdot f\cdot g$, where $c\in K$ is a constant, $f$ is a product of linear polynomials in~$x$ and~$y$ with integer coefficients and $g$
has only finite roots in~$\bN^2$.
\end{theorem}
\begin{proof}
Let~$p = p_1 \cdots p_r$ with $p_i$ irreducible over~$K$.  Assume that $p_1, \ldots, p_m$
have only finitely many zeros in~$\bN^2$ and that $p_i$ with~$i>m$ has infinitely many roots in~$\bN^2$.
Then let~$g = p_1\cdots p_m$. We show that $p_{m+1}, \ldots, p_{r}$ are, up to scalar multiplication, polynomials of the form~$ax+by+c$ with~$a, b, c\in \bZ$.
By Proposition~\ref{LEM:semilinear}, the set $E$ of all nonnegative points $(n, m)$ on the curve~$p(x, y)=0$ is semilinear.
Now suppose that $E$ is infinite. Then if the subset $V_i$ in \eqref{EQ:semilinear} is not contained in a line in~$\bZ^2$ through the origin, then the set $$b_i + \sum_{v\in V_i} v\cdot \mathbb{N}$$ is
Zariski dense in the plane, which is impossible since $E$ is contained in the zero set of a nonzero polynomial.  Thus we see that after refining our decomposition of~$E$ if necessary,
we may assume that each $|V_i| = 1$ for~$i>0$. Let~$q$ be any irreducible factor of~$p$ having infinitely many
zeros in~$\bN^2$. Then there is some $V_i = \{v\}\subseteq \bN^2$ with~$i>0$, such that $q(b_i + vn)=0$
for infinitely many~$n\in \bN$. Write $b_i = (c, d)$ and~$v = (a, b)$. Then $q(c+an, d+bn)=0$ for infinitely many $n\in \bN$ and so $q(c+at, d+bt)=0$ for all $t\in K$.
Hence the linear polynomial $ay-bx-(da-cb)$ divides $q$. Since $q$ is irreducible over~$K$, then $q= \lambda (ay-bx-(da-cb))$
for some constant $\lambda \in K$. This completes the proof.
\end{proof}
The theorem as above cannot be extended to the case when~$d>2$ as shown in the following example.

\begin{example}
Let~$p = x-y + 2z^2 + zy^2$. We claim that $E := \{(n, n, 0)\mid n\in \bN\}$ is the set of
all zeros of $p$ in~$\bN^3$. Suppose that $(a, b, c)$ is another $\bN^3$-point with~$c$ nonzero.
Then $a+2c^2 + cb^2=b$ and so $c(2c+b) = 2c^2+cb^2\le b$ since $a$ is nonnegative. But if $c$ is strictly positive then we must have $2c+b\le c(2c+b)\le b$, which gives $c\le 0$, a contradiction.

Now the corresponding generating function
is equal to $1/((1-x)(1-y))$ which is rational,
but the polynomial $p$ is not of the integer-linear form up to scalar multiplication.
\end{example}

As in the first question, we can show that it is undecidable to test whether
the generating function $F_V$ for an arbitrary algebraic variety~$V$ is D-finite or not.
Let~$P\in \bQ[x_1, \ldots, x_d]$ be any polynomials over~$\bQ$ in $x_1, \ldots, x_d$
and let $V$ be the algebraic variety defined by
\[V := \{(a_1, \ldots, a_d, b, c)\in \overline{\bQ} \mid P(a_1, \ldots, a_d)^2 + (b-c^2)^2=0\}.\]
The undecidability follows from the equivalence that the generating function $F_V$ is D-finite if and only if $P$ has no root in~$\bN^d$.
Clearly, $F_V = 0$ if $P$ has no root in~$\bN^d$ and then it is D-finite.
Now suppose that $P$ has at least one root in~$\bN^d$. Then the generating function $F_V$
is of the form
\[F_V = \sum_{\substack{(n_1, \ldots, n_d, m)\in \bN^{d+1}\\ P(n_1, \ldots, n_d)=0}} x_1^{n_1}\cdots x_d^{n_d} y^{m^2}z^m.\]
It is sufficient to show that $G_V(x_1, \ldots, x_d, y):=F_V(x_1, \ldots, x_d, y, 1)$ is not D-finite.
Clearly, the set
\[\{m\mid \text{$\exists (n_1, \ldots, n_d) \in \bN^d$ such that $g(n_1, \ldots, n_d, m)\neq 0$} \}\]
is the set of square numbers, which is neither finite nor syndetic. Thus $G_V$ is not D-finite by Lemma~\ref{LEM:syndetic}.

\begin{example} Let~$p= x^2-y \in K[x, y]$. Then the associated generating function is~$F(x, y) = \sum_{m\geq 0} x^{m}y^{m^2}$.
Since $p$ is not of the integer-linear form, $F(x, y)$ is not D-finite by Theorem~\ref{THM:d=2}.
Actually, we can show that $F(x, y)$ is transcendentally transcendental. Suppose that $F(x, y)$ is
differentially algebraic. Then it satisfies a nontrivial algebraic differential equation~$Q(x, y, F, D_x(F), \ldots, D_x^r(F)) = 0$, where~$r\in \bN$ and~$Q\in K[z_1, z_2, \ldots, z_{r+3}]$. Note that the evaluation of a power series at $y=2$ gives a ring homomorphism $e_{2}:K[[x, y]]\to K[[x]]$ and we have a commuting square
\[ \begin{array}{ccc} K[[x,y]] & \stackrel{e_{2}}{\longrightarrow} & K[[x]] \\ \big \downarrow & ~ & \big\downarrow\\
K[[x,y]] & \stackrel{e_{2}}{\longrightarrow} & K[[x]],\end{array}\\
 \]
 where both vertical maps are differentiation with respect to $x$.
It follows that $F(x, 2) = \sum_{m\geq 0} 2^{m^2}x^{m}$ is also differentially algebraic.
This leads to a contradiction with the fact proved by Mahler in~\cite[p. 200, Theorem 16]{Mahler1976} on the rate of coefficient growth of a differentially algebraic power series, since $2^{m^2}\gg (m!)^c$ for any positive constant $c$.
\end{example}

This example motivates us to formulate the following conjecture, which can be viewed as an analogue
of the P\'olya-Carlson theorem in the context of algebraic geometry and differential algebra.
\begin{conjecture}
Let $V\subseteq \mathbb{A}^d_K$ be an affine variety over an algebraically closed field $K$ of characteristic zero.  Then
the power series
\[ F_V(x_1, \ldots, x_d) := \sum_{(n_1,\ldots ,n_d)\in V\cap \mathbb{N}^d} x_1^{n_1}\cdots x_d^{n_d}\]
is either rational or transcendentally transcendental.
\end{conjecture}


\bibliographystyle{amsplain}
\def\cdprime{$''$}
\providecommand{\bysame}{\leavevmode\hbox to3em{\hrulefill}\thinspace}
\providecommand{\MR}{\relax\ifhmode\unskip\space\fi MR }
\providecommand{\MRhref}[2]{%
  \href{http://www.ams.org/mathscinet-getitem?mr=#1}{#2}
}
\providecommand{\href}[2]{#2}

%
%
%
%

\end{document}